\documentclass[final,1p,times]{elsarticle}

\usepackage{amsmath,amssymb,amsfonts,amsthm,enumerate,multirow,mathtools}
\usepackage{color}
\usepackage{siunitx}
\usepackage[linesnumbered,ruled]{algorithm2e} 
\usepackage{diagbox} 
\usepackage{cleveref}
\usepackage{placeins} 
\usepackage[table]{xcolor}
\usepackage{graphicx}
\usepackage{tikz}
\usepackage{makecell}

\newtheorem{theorem}{Theorem}
\newtheorem{remark}{Remark}
\newtheorem{conjecture}{Conjecture}
\newtheorem{corollary}{Corollary}
\newtheorem{question}{Question}

\usepackage{lineno}

\journal{TBC}

\begin{document}

\begin{frontmatter}

\title{Uniqueness of solutions in high-energy x-ray based `eigenstrain tomography' and other inverse eigenstrain problems: Counter examples and necessary conditions for well-posedness}

\author[Newc]{CM Wensrich\corref{cor}}
\ead{christopher.wensrich@newcastle.edu.au}
\author[Manc]{S Holman}
\author[Manc]{WBL Lionheart}
\author[Sant]{M Courdurier}
\author[Newc]{RR Jackson}

\cortext[cor]{Corresponding author}

\address[Newc]{School of Engineering, The University of Newcastle, University Drive, Callaghan, NSW 2308, Australia}
\address[Manc]{Department of Mathematics, Alan Turing Building, University of Manchester, Oxford Rd, Manchester M139PL, UK}
\address[Sant]{Department of Mathematics, Pontificia Universidad Católica de Chile, Avda. Vicuña Mackenna 4860, Macul, Santiago, Chile}

\begin{abstract}
Eigenstrain tomography combines diffraction-based strain measurement with elasticity theory to reconstruct full three-dimensional residual stress fields within solids. Notwithstanding a number of recent examples, the uniqueness of such reconstructions has not yet been clearly established. In this paper, we examine the underlying inverse problem in detail and construct explicit counterexamples demonstrating non-uniqueness for a recent implementation of x-ray eigenstrain tomography involving reconstruction from a single measured component of strain.  We follow on to explore minimum conditions for well-posedness and conclude that the full elastic strain tensor within an isotropic sample can be uniquely reconstructed from three measured components; specifically the three shear components, or the three diagonal components. We further prove two key results related to eigenstrain reconstruction in a general sense; 1. That any possible residual stress field can be generated by a diagonal eigenstrain and 2. That residual stress fields exist that cannot be generated by isotropic eigenstrains.  Together, these findings establish rigorous minimum experimental and computational requirements for well-posed eigenstrain tomography techniques and inverse eigenstrain problems in general.
\end{abstract}


\end{frontmatter}

\section{Introduction}

Developments in energy-resolved neutron transmission and high-energy synchrotron x-ray diffraction techniques have given birth to a range of tomographic approaches for imaging elastic strain within polycrystalline solids (e.g. \cite{uzun_tomographic_2024,uzun_voxel-based_2023,gregg_tomographic_2018,hendriks_bragg-edge_2017,hendriks_tomographic_2019,abbey_neutron_2012,korsunsky_strain_2011,korsunsky_principle_2006,wensrich_bragg-edge_2016,wensrich_direct_2024,wensrich_general_2025,uzun_critical_2025}).  While these approaches vary, the unifying feature is the common aim of reconstructing the full elastic strain tensor (and hence stress) at every point within a sample from a set of diffraction-based measurements.

\begin{figure}
    \centering
    \includegraphics[width=0.7\linewidth]{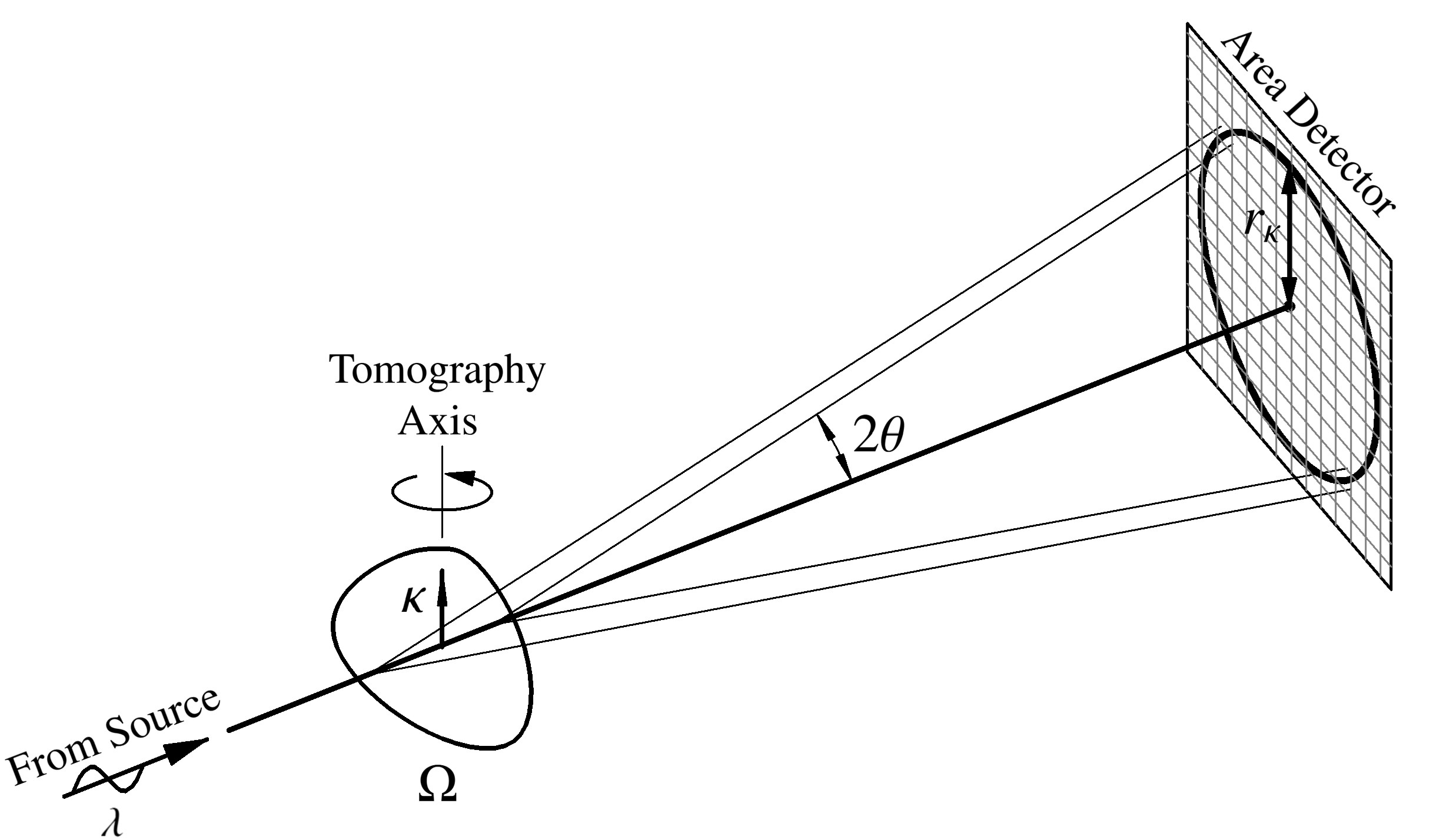}
    \caption{At high energies, x-ray diffraction angles become small and the radius of Debye Scherrer rings can be used to measure elastic strain in the (approximately) transverse direction relative to the incident beam.  In the vertical direction $\kappa$, the only component that contributes to this measurement is invariant to the rotation angle and can be reconstructed through standard scalar filtered back projection \cite{korsunsky_strain_2011,uzun_critical_2025}.}
    \label{fig:XrayGeometry}
\end{figure}

Of particular interest to this paper is a technique involving high-energy synchrotron x-ray diffraction described by Korsunsky \emph{et al.} \cite{korsunsky_strain_2011,korsunsky_principle_2006}.  Also described conceptually by Hendriks \emph{et al.} \cite{hendriks_bayesian_2020} the technique involves the measurement of elastic strain within a sample through relative shift in Debye Scherrer rings formed by the diffraction of high-energy x-rays (see Figure \ref{fig:XrayGeometry}).  The geometry of these rings is defined by Bragg's law
\[
m\lambda=2d\sin\theta
\]
where $m\in\mathbb{N}$ and $\theta$ is the diffraction angle formed by x-rays of wavelength $\lambda$ scattering from lattice planes with spacing $d$.  In a sample subject to residual stress, lattice spacing is related to elastic strain and it follows that
\[
\epsilon_{\kappa\kappa}=\frac{d_\kappa-d_0}{d_0}\approx\frac{\theta_0-\theta_\kappa}{\theta_\kappa}\approx\frac{r_0-r_\kappa}{r_\kappa}
\]
where $\epsilon_{\kappa\kappa}=\epsilon_{ij}\kappa_i\kappa_j$ is the normal component of elastic strain measured in the $\kappa$ direction, $d_\kappa,\theta_\kappa,r_\kappa$ are the lattice spacing, diffraction angle and radius corresponding to $\kappa$, and $d_0,\theta_0,r_0$ are equivalent values for reference/unstressed material.  Note that for high-energy x-rays, $2\theta$ is small and $\kappa$ is approximately transverse with respect to the incident beam.

A strain measurement made in this way corresponds to an average along the entire ray-path through the sample and the collection of many rays of this form naturally results in a tomography problem for elastic strain based on the \emph{Transverse Ray Transform} \cite{sharafutdinov_integral_2012,lionheart_diffraction_2015}. Note that this tomography problem involves the reconstruction of a rank-2 tensor field at every point within a sample, however there is a key simplifying feature as follows \cite{lionheart_diffraction_2015};

For a given tomographic axis of rotation, the measured component of normal elastic strain along the axis is independent of the rotation angle.  With reference to Figure \ref{fig:XrayGeometry}, this means that for the rotation axis shown, it is possible to tomographically reconstruct the $\epsilon_{\kappa\kappa}$ through traditional parallel-beam (scalar) filtered back projection.  If we repeat this process \emph{six} times with sufficiently independent rotation axes we can solve for the six unique components of the elastic strain tensor at every point.\footnote{The required condition is that not all rotation axes can lie on the surface of any given projective conic (e.g. cone) \cite{desai_explicit_2016,lionheart_diffraction_2015,wensrich_force_2014}}

While this is a simple concept, it defines a formidable experiment involving six separate tomographic scans that require careful relative alignment. In this setting, \emph{eigenstrain tomography} promises a simplifying path.  As the name would suggest, eigenstrain tomography relies heavily on the eigenstrain approach in elasticity theory in the reconstruction process.  In essence, this approach provides a path to impose physical constraints such as mechanical equilibrium and known boundary conditions on the reconstruction process.  The hope is that these constraints can reduce or minimise the number of components that need to be measured, with the others calculated through elasticity theory.  The central question of this paper is; ``What are the minimum requirements of this process?"  In other words, how many components do we need to measure to ensure the final reconstruction is robust?

We begin with an overview of eigenstrain theory as it relates to this problem.

\section{Eigenstrain theory of residual stress}

Originally used to model and understand the effect of inclusions and defects within elastic bodies (e.g. \cite{lin_elastic_1973,mura_elastic_1976,eshelby_determination_1957,kinoshita_elastic_1971}), eigenstrain has become a central concept in our understanding of residual stress in solids \cite{korsunsky_eigenstrain_2008,korsunsky_teaching_2017,hill_determination_1996,wensrich_well-posedness_2025}.  A summary of the approach is as follows;

Consider a finite deformable solid body represented by the compact domain $\Omega\in\mathbb{R}^3$ with unit normal $n$ on a Lipchitz boundary $\partial\Omega$.  This body is subject to a continuous infinitesimal displacement field $X\in C^2(\mathcal{S}^1,\Omega)$, from which we can compute the total strain at any point in $\Omega$ as
\begin{equation} \label{eq:Decomp}
\epsilon^T=\nabla_sX=\epsilon+\epsilon^*,
\end{equation}
where $\nabla_s$ is the symmetric gradient operator, $\epsilon\in C^2(\mathcal{S}^2,\Omega)$ is the elastic component of strain and $\epsilon^*\in C^2(\mathcal{S}^2,\Omega)$ is an inelastic eigenstrain.\footnote{See \ref{Sec:Apndx} for notation and spaces used throughout this paper} 

Elastic strain is related directly to stress within the body through Hooke's law
\[
\sigma=C:\epsilon
\]
where $C\in\mathcal{S}^4$ is the usual positive definite tensor of elastic constants.  In the simplest case of a homogeneous isotropic body, all 81 components of $C$ can be expressed in terms of 2 scalar elastic properties (e.g. Young's modulus $E$ and Poisson's ratio $\nu$).  We will assume this is the case throughout the paper.

The eigenstrain component represents various forms of inelastic strain that may be present within the sample (e.g. due to plasticity, phase change, thermal expansion, interference, etc).  This decomposition of strain into elastic and eigen components is the central feature of the eigenstrain approach.

In the case of residual stress we assume no surface traction or body forces are present and mechanical equilibrium is satisfied provided
\begin{align}\label{eq:Equilibrium}
\mathrm{Div}(\sigma)&=0 \text{ on } \Omega, \text{ and}\\
\sigma\cdot n &=  0 \text{ on } \partial\Omega.
\end{align}
It follows from \eqref{eq:Decomp} that the displacement field resulting from a given eigenstrain can be determined as the unique solution of the elliptic boundary value problem
\[
\begin{split}
\mathrm{Div}(C:\nabla_sX)&=\mathrm{Div}(C:\epsilon^*) \quad \text{on} \quad \Omega,\\
(C:\nabla_sX)\cdot n&=(C:\epsilon^*)\cdot n \quad \text{on} \quad \partial\Omega.
\end{split}
\]
Practically, numerical solutions for $X$ can be found for a given $\epsilon^*$ using established finite-element techniques from which residual stress within the body can be determined as
\[
\sigma=C:(\nabla_sX-\epsilon^*).
\]

Armed with this `forward' calculation (i.e. the mapping from $\epsilon^*$ to $\sigma$), the eigenstrain tomography problem can be generally stated as follows;

\begin{question}
    Say we have incomplete data related to the stress field within $\Omega$.  For example, we may have knowledge of one or more components of $\epsilon$ or $\sigma$ at every point in $\Omega$, or knowledge of $X$ on several surfaces within $\Omega$.  Can we solve for a corresponding $\epsilon^*$ that uniquely defines $\sigma$ at every point within $\Omega$?
\end{question}

Philosophically this idea extends beyond diffraction-based measurements; any measurement that defines a suitable inverse problem for eigenstrain could be coined eigenstrain tomography, however generically the phrase `eigenstrain approach' is more commonly used in other circumstances (e.g. \cite{hill_determination_1996, korsunsky_eigenstrain_2008,korsunsky_teaching_2017}).

The central issue with an inverse problem such as this is the uniqueness of solutions; the less data we have, the more we should expect ambiguity in resulting solutions.  It is common to approach a problem such as this as a numerical optimisation and from this perspective the uniqueness issue can be masked by inherent regularisation within the algorithm employed.  For example, when faced with ambiguity, least-squares minimisation based on the Moore-Penrose pseudo-inverse will determine the solution with the minimum $L^2$ norm.  While this solution may be physically feasible, there is no reason to suspect that it is true; other regularisation schemes will arrive at other solutions that are just as feasible.

From this perspective it is easy to overlook the uniqueness issue.  Uzun \emph{et al.} \cite{uzun_tomographic_2024} provides a prominent recent example that we discuss in detail in the following section.

\section{Ill-posedness of 1 and 2 component eigenstrain tomography}

Using a high-energy x-ray diffraction technique, Uzun \emph{et al.} \cite{uzun_tomographic_2024} examined residual stress within a rectangular additively manufactured nickel-based super-alloy sample. This involved eigenstrain tomography in the form of a full-field reconstruction of the stress tensor within their sample from knowledge of single component of elastic strain measured using the technique described earlier.

Without loss of generality, if we assume this known component is $\epsilon_{33}$ we can demonstrate the non-uniqueness of their solution by constructing a set of counter examples in the form of non-zero strain fields with $\epsilon_{33}=0$ at every point in a given domain. Hypothetically, any linear combination of these counter-examples can be added to their solution without altering their measured component.

We can actually do more than this; 

\begin{remark}{It is possible to construct otherwise non-zero strain fields with both $\epsilon_{22}=\epsilon_{33}=0$ throughout a given domain.}
\end{remark}

\noindent
We demonstrate this as follows;

For simplicity, let $\Omega$ be the cube $[0,2\pi]^3$.  The zero-traction boundary condition on $\Omega$ amounts to $\sigma_{ij} = 0$ for $x_j = 0$ or $2\pi$. 

Within $\Omega$, we can construct an appropriate divergence-free stress $\sigma$ using a Maxwell potential of the form
\begin{equation} \label{eq:Maxwell}
\sigma=\nabla \times \left(\nabla \times \mathrm{diag}(\Lambda)\right)^T
        =\begin{bmatrix}
        \partial_3^2\Lambda_2+\partial_2^2\Lambda_3& -\partial_1\partial_2\Lambda_3 & -\partial_1\partial_3\Lambda_2 \\
            -\partial_1\partial_2\Lambda_3 &  \partial_1^2\Lambda_3+\partial_3^2\Lambda_1 & -\partial_2\partial_3\Lambda_1 \\
            -\partial_1\partial_3\Lambda_2 & -\partial_2\partial_3\Lambda_1 &  \partial_2^2\Lambda_1+\partial_1^2\Lambda_2
        \end{bmatrix},
\end{equation}
where $\Lambda\in C^4(\mathcal{S}^1,\Omega)$ is an arbitrary vector potential and the curl-curl operator acts on both the rows and columns of the diagonal tensor. Note that any solution to \eqref{eq:Equilibrium} can be written in this form \cite{langhaar_three-dimensional_1954}. Since the material is isotropic, Hooke's law gives the components $\epsilon_{33}$ and $\epsilon_{22}$ as
\[
\begin{split}
\epsilon_{22} & = \frac{1}{E} (\sigma_{22} - \nu (\sigma_{11}+\sigma_{33})),\\
\epsilon_{33} & = \frac{1}{E} (\sigma_{33} - \nu (\sigma_{11}+\sigma_{22})).
\end{split}
\]
Setting both of these equal to zero and substituting from \eqref{eq:Maxwell} gives a pair of partial differential equations for $\Lambda$, and we are interested in finding solutions subject to the zero-traction boundary conditions. The system can be recast in a slightly simplified way if we replace $\epsilon_{22}=0$ by the expression $\sigma_{33}-\sigma_{22}=0$ and obtain
\[
\begin{split}
\sigma_{33}-\nu(\sigma_{11} + \sigma_{22}) & = 0,\\
\sigma_{33}-\sigma_{22} & = 0.
\end{split}
\]
In terms of $\Lambda$ this becomes
\begin{equation}\label{eq:epsys}
\begin{split}
\partial_1^2 (\Lambda_2 - \nu \Lambda_3) + \partial_2^2(\Lambda_1 - \nu \Lambda_3) - \nu \partial_3^2 (\Lambda_1 + \Lambda_2) & = 0,\\
\partial_1^2(\Lambda_2-\Lambda_3) + \partial_2^2 \Lambda_1 - \partial_3^2 \Lambda_1 & = 0.
\end{split}
\end{equation}
Let us pose an ansatz for $\Lambda$ in the form of cosine series
\begin{equation}\label{eq:Ansatz}
\Lambda_i =  \sum_{j,k,l=1}^\infty a^i_{jkl} \cos(j x_1)\cos(k x_2) \cos(l x_3).
\end{equation}
As per the usual approaches in Fourier analysis, the terms in the sum are orthogonal and it can be shown that the map from the set of possible coefficients $a^i_{jkl}$ to $\sigma$ is injective (this is not the case if the sum starts at zero).  Substituting \eqref{eq:Ansatz} into \eqref{eq:epsys} we obtain
\[
\begin{split}
(k^2 - \nu l^2) a^{1}_{jkl}+(j^2-\nu l^2) a^2_{jkl} - \nu (j^2 + k^2) a^3_{jkl} & = 0\\
(k^2 - l^2) a^1_{jkl} + j^2 a^2_{jkl} - j^2 a^3_{jkl} & = 0.
\end{split}
\]
These equations are equivalent to
\begin{equation}\label{eq:bcform}
\begin{bmatrix}
    a^1_{jkl} \\ a^2_{jkl} \\ a^3_{jkl}
\end{bmatrix}
=
b_{jkl}
\begin{bmatrix}
    k^2 - \nu l^2\\
    j^2 - \nu l^2\\
    -\nu (j^2 + k^2)
\end{bmatrix}
\times
\begin{bmatrix}
    k^2 - l^2\\
    j^2\\
    -j^2
\end{bmatrix}
=
b_{jkl}
\begin{bmatrix}
    j^2(-(1-\nu)j^2 + \nu k^2+\nu l^2)\\
    k^2((1-\nu) j^2 - \nu k^2 +\nu l^2)\\
    l^2((1-\nu) j^2 + \nu k^2-\nu l^2)
\end{bmatrix}
\end{equation}
for some coefficients $b_{jkl}$. In summary, whenever $a^i_{jkl}$ are given by \eqref{eq:bcform} for arbitrary $b_{jkl}$ we have $\epsilon_{22} = \epsilon_{33} = 0$.\footnote{As an alternate approach, it is possible to construct arbitrary Maxwell potentials of the form required.  With inspiration from \eqref{eq:bcform}, suppose we define a differential operator $\mathcal{A}$ by
\[
\mathcal{A} = \begin{bmatrix}
    \partial_2^2 - \nu\partial_2^2\\
    \partial_1^2 - \nu \partial_3^2\\
    -\nu \partial_1^2 - \nu \partial_2^2
\end{bmatrix}
\times
\begin{bmatrix}
    \partial_2^2 - \partial_3^2\\
    \partial_1^2\\
    -\partial_1^2
\end{bmatrix}
=
\begin{bmatrix}
    \big(-(1-\nu)\partial_1^2 + \nu \partial_2^2 + \nu \partial_3^2\big) \partial_1^2\\
    \big((1-\nu)\partial_1^2 - \nu \partial_2^2 + \nu \partial_3^2\big)\partial_2^2\\
    \big((1-\nu)\partial_1^2 + \nu \partial_2^2 - \nu \partial_3^2\big) \partial_3^2
\end{bmatrix}.
\]
Then for any $a \in C^6(\mathbb{R}^3)$, $\Lambda = \mathcal{A}a$ automatically satisfies \eqref{eq:epsys}. Thus we can construct strain fields with arbitrary support satisfying $\epsilon_{22} = \epsilon_{33} = 0$.}

Now let us consider the boundary conditions. Note that, due to the form of the ansatz, on $\partial\Omega$ the terms in the boundary condition that involve mixed derivatives are already zero by construction. Thus the boundary conditions reduce to a single equation on each component of the boundary involving the diagonal components of $\sigma$. The conditions on opposite sides of the cube are the same by periodicity and so there are only three equations;
\[
\sum_{j=1}^\infty l^2 a^{2}_{jkl} + k^2 a^{3}_{jkl} = 0, \ \sum_{k=1}^\infty j^2 a^{3}_{jkl} + l^2 a^{1}_{jkl} = 0, \ \sum_{l=1}^\infty k^2 a^{1}_{jkl} + j^2 a^{2}_{jkl} = 0.
\]
Using \eqref{eq:bcform}, these three boundary conditions become equations for $b_{jkl}$ which simplify significantly. Indeed, the first equation becomes
\[
\begin{split}
0 & = \sum_{j=1}^\infty \Big ( l^2 k^2((1-\nu) j^2 - \nu k^2 +\nu l^2) + k^2  l^2((1-\nu) j^2 + \nu k^2-\nu l^2) \Big ) b_{jkl}\\
& = 2\sum_{j=1}^\infty j^2 k^2 l^2 (1-\nu) b_{jkl}
\end{split}
\]
which is equivalent to
\begin{equation} \label{eq:j_eqn}
    \sum_{j=1}^\infty j^2 b_{jkl} = 0,
\end{equation}
for all $k,l$.
Similar calculations show
\begin{equation}\label{eq:kl_eqn}
    \sum_{k=1}^\infty k^2 b_{jkl} = 0, \ \sum_{l = 1}^\infty l^2 b_{jkl} = 0,
\end{equation}
for all $j,l$ and $j,k$ respectively.

Supposing that we restrict to $j,k,l \leq N$ for some $N\in\mathbb{N}$, these amount to $3N^2$ equations. Since there are $N^3$ corresponding coefficients $b_{jkl}$, we see that for $N>3$ we are guaranteed nontrivial solutions.

To implement this scheme, we begin by flattening $b_{jkl}$ into a one-dimensional array with $N^3$ elements through an appropriate indexing function.  Our choice maps element $b_{jkl}$ with $j,k,l\in[1,N]$ to $B_i$ with $i\in[1,N^3]$ where
\[
i=j+(k-1)N+(l-1)N^2.
\]
The inverse of this mapping is
\[
\begin{split}
    j&=\big[(i-1) \text{ mod } N\big]+1,\\
    k&=\big[\lfloor \tfrac{i-1}{N} \rfloor \text{ mod } N\big]+1,\\
    l&=\lfloor \tfrac{i-1}{N^2} \rfloor +1.
\end{split}
\]

We can then assemble \eqref{eq:j_eqn} and \eqref{eq:kl_eqn} into a system of the form
\[
A_{\big[3N^2\times N^3\big]}B_{\big[N^3\times1\big]}=0_{\big[3N^2\times1\big]},
\]
and seek non-trivial solutions for $B$.

\begin{figure} [!ht]
    \centering
    \includegraphics[height=0.18\linewidth]{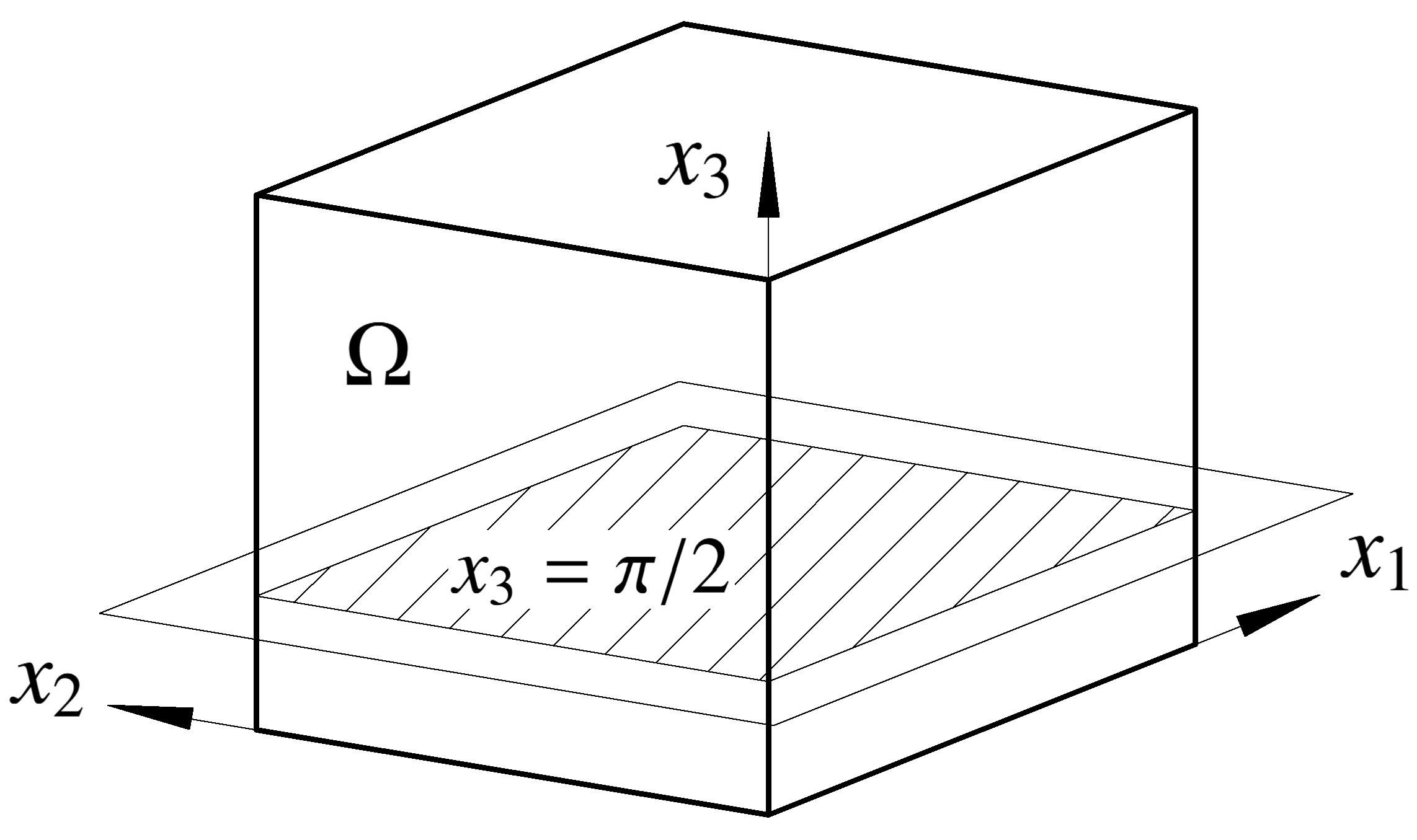}\\
    \includegraphics[width=0.45\linewidth]{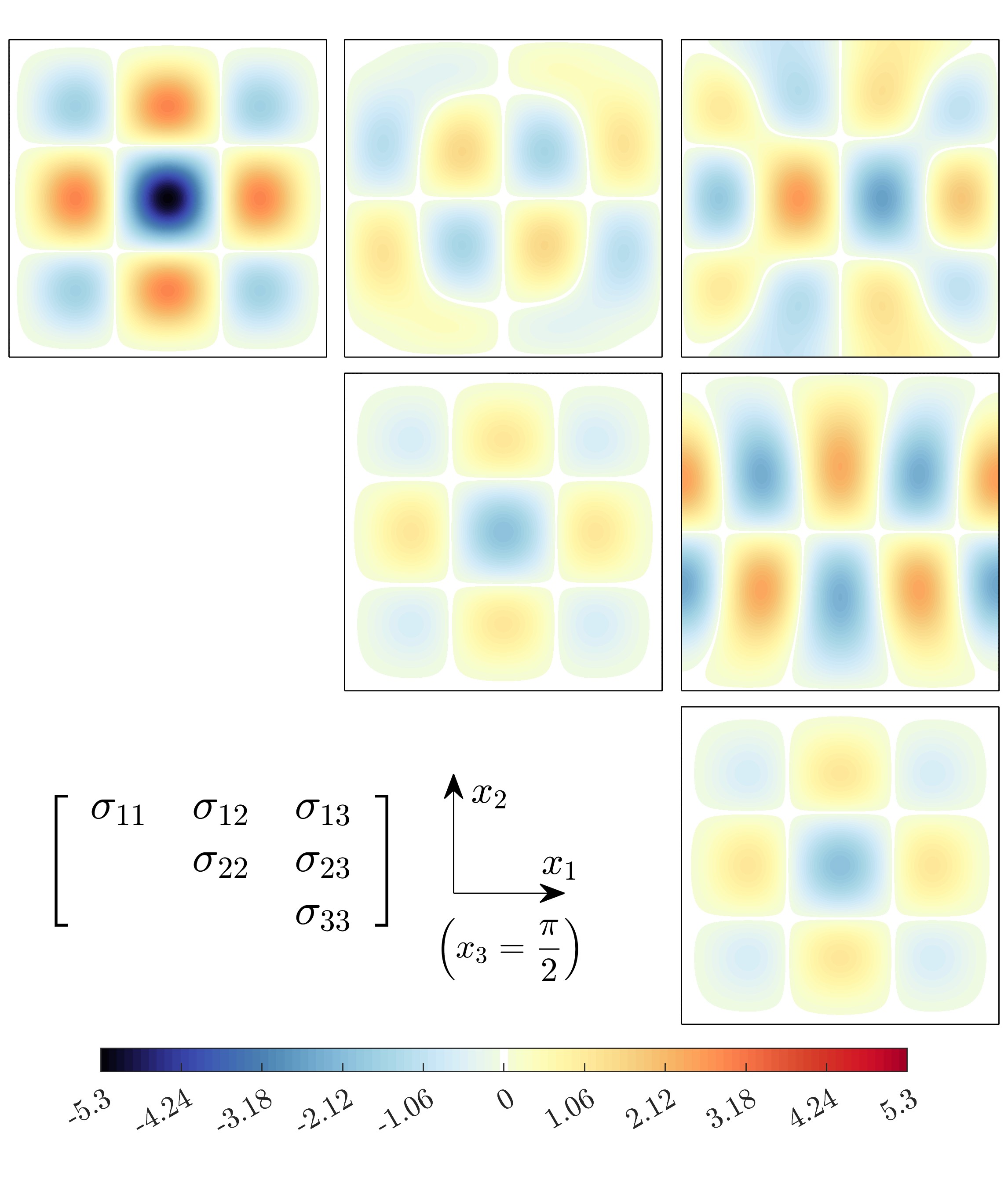}
    \includegraphics[width=0.45\linewidth]{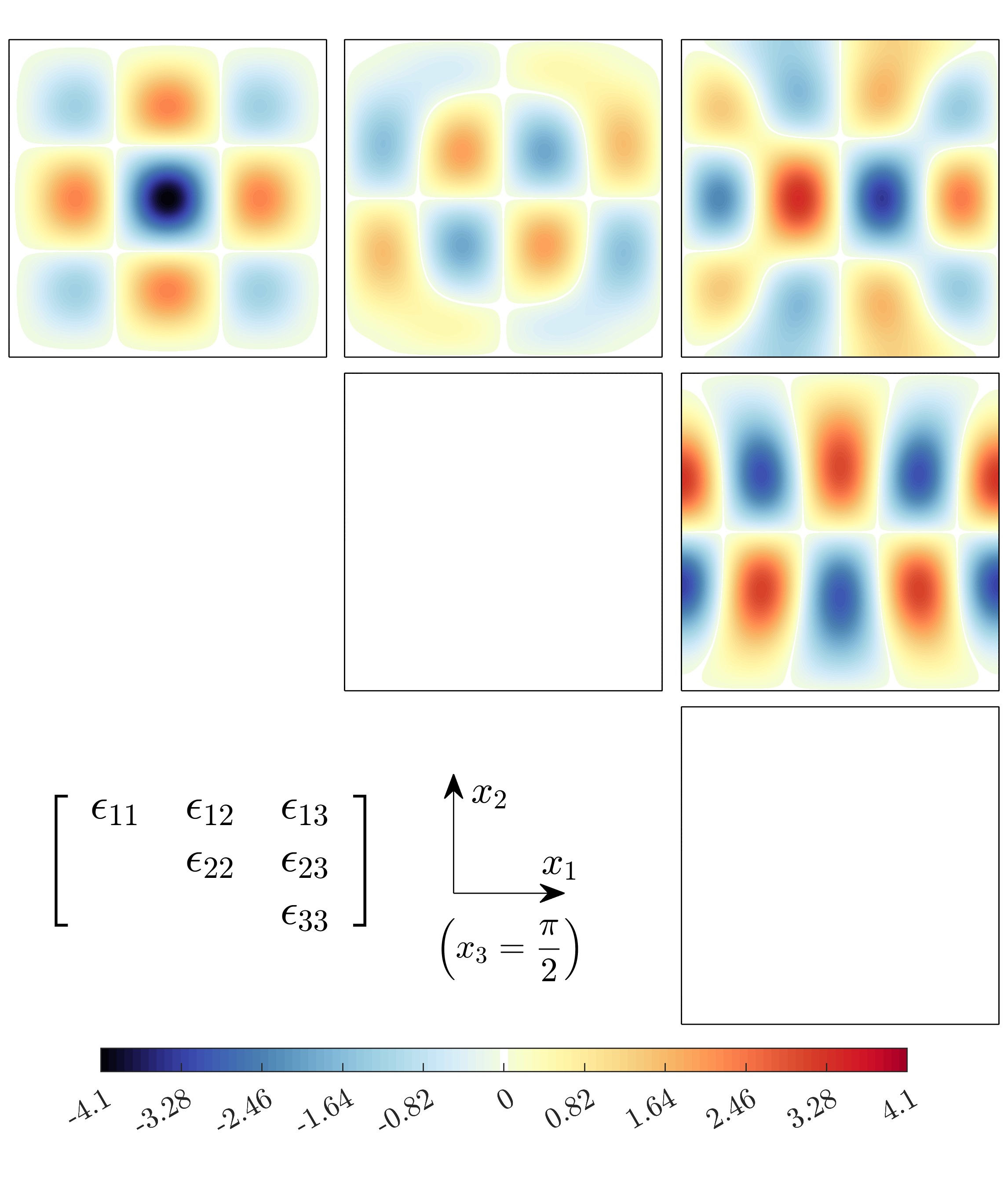}
    \caption{Stress and the corresponding elastic strain within the cube $[0,2\pi]^3$ over a cross section at $x_3=\pi/2$ (with $E=1$ and $\nu=0.28$).  This corresponds to the single null vector of $A$ in the case of $N=2$.  It represents a simple example of a residual stress field on the cube with $\epsilon_{22}=\epsilon_{33}=0$ while simultaneously satisfying equilibrium and a zero-traction boundary condition.}
    \label{fig:firstbasisvector}
\end{figure}

With $E=1$ and $\nu=0.28$, this process was carried out in MATLAB using the intrinsic function `\texttt{null}' to generate an orthonormal basis for the null space of $A$.  A single non-trivial solution was actually found for $N=2$ and its corresponding stress and strain field over a planar slice at $x_3=\pi/2$ is shown in Figure \ref{fig:firstbasisvector}.  This represents the simplest non-trivial solution of this form that can be found for this domain.

\begin{figure} [!t]
    \centering
    \includegraphics[height=0.15\linewidth]{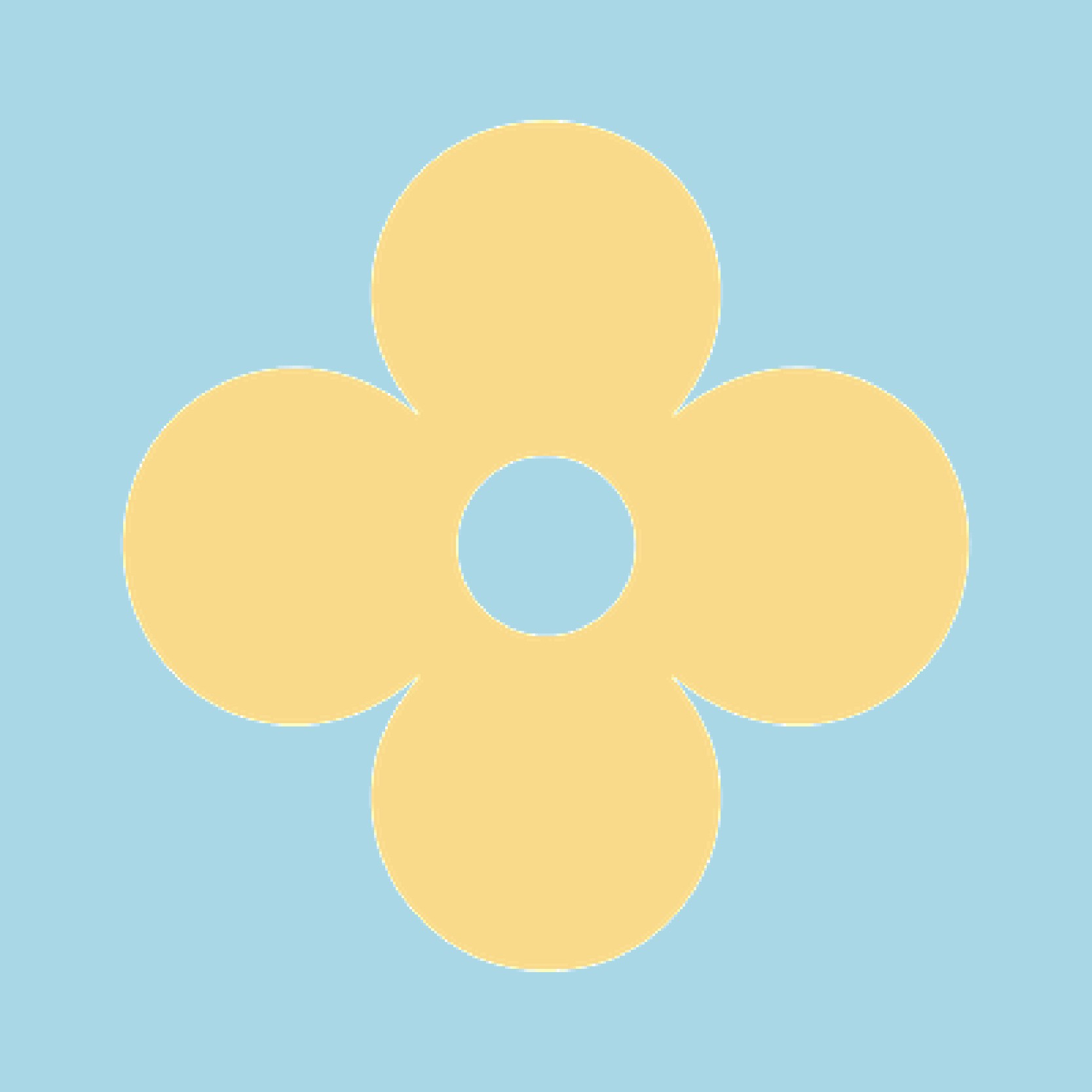}
    \put(-135,45){\small{Target Image ($\sigma_{11}$):}}\\
    \includegraphics[width=0.45\linewidth]{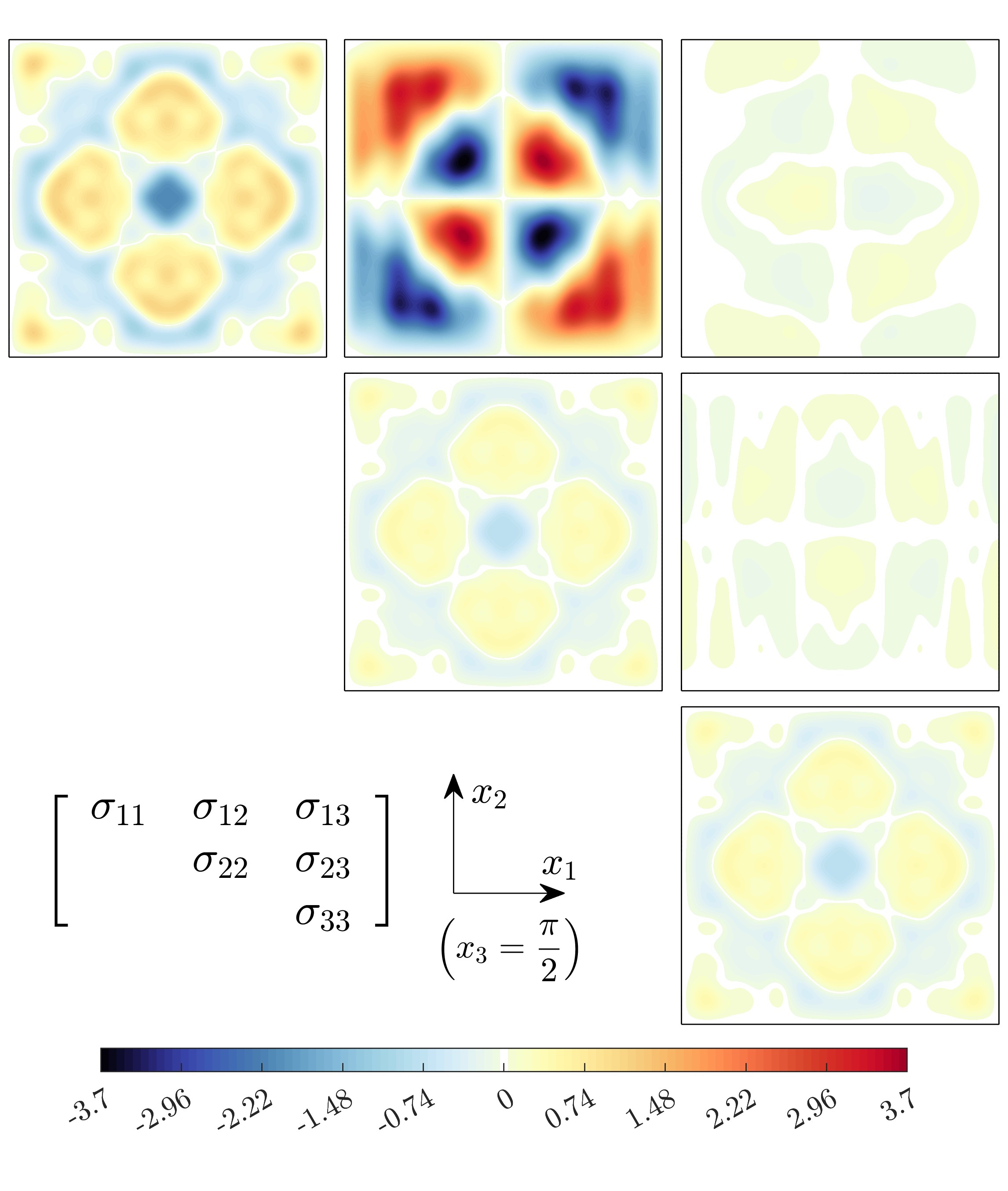}
    \includegraphics[width=0.45\linewidth]{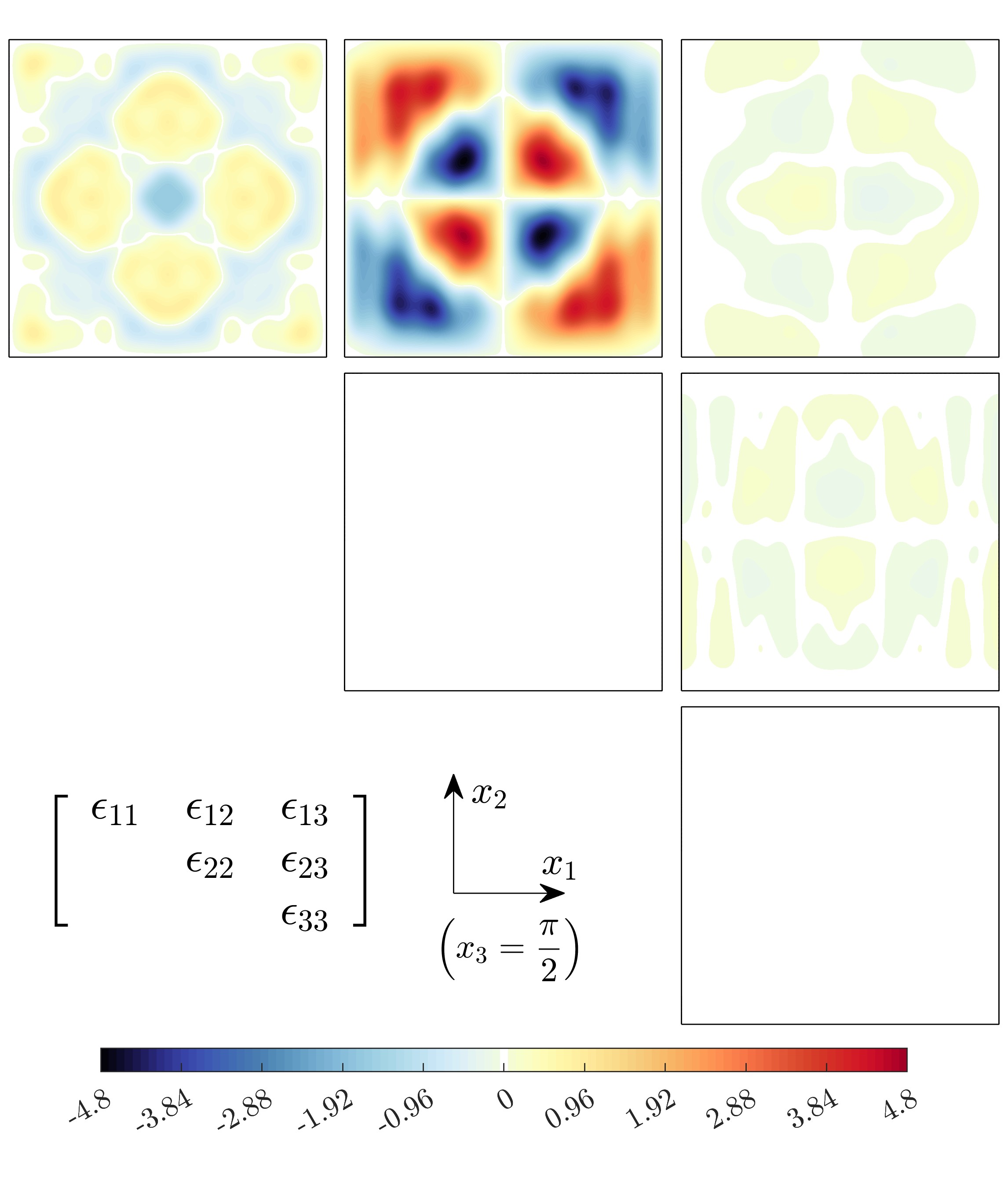}
    \caption{Similar to Figure \ref{fig:firstbasisvector}, a residual stress/strain field shown on the planar slice at $x_3=\pi/2$, but now corresponding to a weighted sum of the 343 null basis fields for $N=8$.  Coefficients in the sum have been chosen through a least-squares process fitting $\sigma_{11}$ to the binary ($\pm1$) image shown.  Note that all boundary conditions are satisfied and $\epsilon_{22}=\epsilon_{33}=0$ at every point in this cube.}
    \label{fig:FittedField}
\end{figure}

For larger $N$, the total number of null-basis elements was $(N-1)^3$. To demonstrate the wide variety of possibilities this affords, Figure \ref{fig:FittedField} shows a residual stress field in the span of the null fields for $N=8$ with $\sigma_{11}$ over the same slice as before fitted to a binary ($\pm1$) image of a flower using a simple Moore-Penrose least squares process.

To summarise, the existence of null-strain fields of this form poses a significant issue for the approach demonstrated by Uzun \emph{et al.} \cite{uzun_tomographic_2024}.  With knowledge of only the $\epsilon_{33}$ component, any strain field of the form shown in Figures \ref{fig:firstbasisvector} or \ref{fig:FittedField} (and many others with $\epsilon_{22}\ne0$) is unobservable by their method.  Conversely, their reconstructed stress field is likely to contain non-physical components of this form implicitly chosen through the regularisation process inherent in their numerical method.

And so, from this perspective we ask the question: \emph{How many components of $\epsilon$ must be measured to uniquely define $\sigma$?}  In this regard, we pose the following conjecture;

\begin{conjecture} \label{Con:Conj1}
    Any elastic strain field, $\epsilon \in C^2(\mathcal{S}^2,\Omega)$, in an isotropic media that satisfies equilibrium and a zero-traction boundary condition over some compact domain $\Omega \subset \mathbb{R}^3$ is uniquely determined from three distinct components.
\end{conjecture}

It is possible that any choice of 3 distinct components is sufficient, however, some are more practical than others.  We consider two versions; 
\begin{enumerate}
    \item{Reconstruction of the diagonal components of $\epsilon$ from the shear components}
    \item{Reconstruction of the shear components of $\epsilon$ from the diagonal components}
\end{enumerate}

\section{Recovery of $\epsilon$ from shear components}

While not necessarily practical from an experimental point of view, we begin by considering the scenario where the three off-diagonal components of $\epsilon$ are known and we seek to solve for the diagonal components of strain (and hence stress).

From Hooke's law in an isotropic material, shear stress components are directly proportional to and defined by the corresponding shear strain components.  So from \eqref{eq:Equilibrium}, and assuming all shear components are zero, we can solve for unique diagonal components of stress provided
\begin{align} \label{eq:DiagonalEquilibrium}
    \partial_1\sigma_{11}&=0 \nonumber \\
    \partial_2\sigma_{22}&=0\\
    \partial_3\sigma_{33}&=0 \nonumber
\end{align}
only has a trivial solution when constrained by the usual zero-traction boundary condition.  This boundary condition on the remaining diagonal components can be written 
\[\sigma_{11}n_1=\sigma_{22}n_2=\sigma_{33}n_3=0
\]
on $\partial \Omega$. 

Noting that \eqref{eq:DiagonalEquilibrium} implies each diagonal component is constant in its respective coordinate direction, it follows that $\sigma=0$ provided the span of the normal vectors is all of $\mathbb{R}^3$ (true for a compact domain).  The argument is as follows;

Consider a point on the boundary with some $n_i\ne0$.  The boundary condition then implies $\sigma_{ii}=0$ (no summation implied) at all points along the line through $\Omega$ with the same $x_j$ and $x_k$ ($j,k\ne i$).  Repeating this process for all boundary points in the three coordinate directions forces $\sigma=0$ throughout $\Omega$.

Hence the three shear components of strain uniquely determine the whole tensor.

\section{Recovery of $\epsilon$ from normal components}

We now turn our attention to the more practical problem of determining the shear components from measurements of the three diagonal components of strain.  From Hooke's law, in an isotropic homogeneous material we can calculate the three diagonal components of stress from the diagonal components of elastic strain.  Assuming these diagonal components are zero, the remaining shear components of stress will then be uniquely defined by \eqref{eq:Equilibrium} provided
\begin{align} \label{eq:ShearEquilibrium}
    \partial_2\sigma_{12}+\partial_3\sigma_{13}&=0 \nonumber\\
    \partial_1\sigma_{12}+\partial_3\sigma_{23}&=0 \\
    \partial_1\sigma_{13}+\partial_2\sigma_{23}&=0 \nonumber
\end{align}
only has a trivial solution when subject to the usual zero-traction boundary condition which can now be written
\begin{equation} \label{eq:BCv1}
        \begin{bmatrix}
        0 & \sigma_{12} & \sigma_{13}\\
        \sigma_{12} & 0 &  \sigma_{23}\\
        \sigma_{13} & \sigma_{23} & 0
    \end{bmatrix}
    \cdot n =0 \text{ on } \partial \Omega.
\end{equation}
The system given in \eqref{eq:ShearEquilibrium} defines interesting structures related to Killing and conformal Killing vector fields that we have dubbed \emph{poorly} vector fields.  See \ref{apdx:poorly} for further discussion.

Writing the unknown vector as $\tau=\begin{bmatrix} \sigma_{23} &\sigma_{13} &\sigma_{12} \end{bmatrix}^T$, the linear operator defined by \eqref{eq:ShearEquilibrium} has principal symbol
\[
\mathcal{L}(\xi) = 
    \begin{bmatrix}
        0  &  \xi_3 & \xi_2 \\
        \xi_3 & 0 & \xi_1\\
        \xi_2 & \xi_1 & 0
    \end{bmatrix},
\]
with determinant $\det( \mathcal{L}) = 2\xi_1\xi_2\xi_3$.  In other words the system is not elliptic; the principal symbol is not invertible on characteristics consisting of the coordinate planes. In this context, we now seek a non-trivial general solution.

\subsection{General solutions for \eqref{eq:ShearEquilibrium}}

Differentiating each of \eqref{eq:ShearEquilibrium} in turn provides
\[
\begin{split}
   \partial_1\partial_2\sigma_{12} &= -\partial_1\partial_3\sigma_{13},\\
    \partial_1\partial_3\sigma_{13} &= -\partial_2\partial_3\sigma_{23},\\
    \partial_2\partial_3\sigma_{23} &= -\partial_1\partial_2\sigma_{12},
\end{split}
\]
which together imply
\begin{equation}\label{eq:Vanishing2ndDerivatives}
\partial_1\partial_2\sigma_{12}=\partial_1\partial_3\sigma_{13}=\partial_2\partial_3\sigma_{23}=0.
\end{equation}
Integrating twice, we can write each shear component in the form
\begin{equation} \label{eq:sigmainitgen}
\begin{split}
    \sigma_{12}&=\phi_{12}(x_1,x_3)+\psi_{12}(x_2,x_3)\\
    \sigma_{13}&=\phi_{13}(x_1,x_2)+\psi_{13}(x_2,x_3)\\
    \sigma_{23}&=\phi_{23}(x_1,x_2)+\psi_{23}(x_1,x_3)
\end{split}
\end{equation}
for arbitrary $\phi_{ij},\psi_{ij}\in C^3(\mathbb{R}^2)$.\footnote{Note that there is a slight abuse of notation here. The functions $\phi_{ij},\psi_{ij}$ are actually defined over the domain $\Omega\subset\mathbb{R}^3$, however they each have a derivative in one of the coordinate directions equal to zero.  Our notation makes clear which direction this is.  The same is true for other functions such as $\phi,\psi$ and $U_i$ which are appear later in the paper. Note also that in \eqref{eq:SeanGeneral} there are integrals with limits from $-\infty$.  A more robust version would integrate in from a corresponding point on $\partial\Omega$, however we opted to abuse the notation rather than increase complexity.}  Substituting back into \eqref{eq:ShearEquilibrium} we obtain
\begin{align}
\label{eq:PhiPsiConstraints}
    \partial_2\psi_{12}+\partial_3\psi_{13}&=0,\nonumber\\
    \partial_1\phi_{12}+\partial_3\psi_{23}&=0,\\
    \partial_1\phi_{13}+\partial_2\phi_{23}&=0.\nonumber
\end{align}

Choosing $\phi_{12}=\phi, \psi_{12}=\psi$ and $\phi_{13}=\omega$ and solving for the remaining functions we arrive at the following form of a general solution;
\begin{align}
\label{eq:SeanGeneral}
\sigma_{12} & = \phi(x_1,x_3) + \psi(x_2,x_3), \nonumber\\
\sigma_{13} & = - \int_{-\infty}^{x_3} \partial_2 \psi(x_2,s) \ \mathrm{d}s + \omega(x_1,x_2),\\
\sigma_{23} & = - \int_{-\infty}^{x_3} \partial_1 \phi(x_1,s) \ \mathrm{d}s- \int_{-\infty}^{x_2} \partial_1 \omega(x_1,s) \ \mathrm{d}s. \nonumber
\end{align}

Alternatively, we can choose
\[
\begin{split}
\phi_{12}&=\partial_3 U_2(x_1,x_3),\\
\psi_{13}&=\partial_2 U_1(x_2,x_3),\\
\phi_{23}&=\partial_1 U_3(x_1,x_2),
\end{split}
\]
for arbitrary $U_i\in C^4(\mathbb{R}^2)$, and from \eqref{eq:PhiPsiConstraints} we then require
\[
\begin{split}
    \psi_{12}&=-\partial_3 U_1(x_2,x_3),\\
    \phi_{13}&=-\partial_2 U_3(x_1,x_2),\\
    \psi_{23}&=-\partial_1 U_2(x_1,x_3),
\end{split}
\]
which provides an alternate form of the general solution
\begin{equation}
    \label{eq:ChrisGeneral}
    \tau=
    \begin{bmatrix}
        \sigma_{23}\\
        \sigma_{13}\\
        \sigma_{12}
    \end{bmatrix}
    =
    \begin{bmatrix}
        \partial_1\left(-U_2(x_1,x_3)+U_3(x_1,x_2)\right)\\
        \partial_2\left(-U_3(x_1,x_2)+U_1(x_2,x_3)\right)\\
        \partial_3\left(-U_1(x_2,x_3)+U_2(x_1,x_3)\right)
    \end{bmatrix}
\end{equation}

These two versions, \eqref{eq:SeanGeneral} and \eqref{eq:ChrisGeneral}, are equivalent parameterisations of the same general solution in terms of three arbitrary functions on the plane.

\subsection{Boundary conditions on an arbitrary convex domain}

Consider $\Omega$ as an arbitrary convex domain with with a $C^1$ boundary. For the time being we restrict to a convex domain for simplicity;  \ref{apdx:DistProof} provides a more general proof for arbitrary compact Lipschitz $\Omega$.  

In terms of $\tau$, the boundary condition can be written
\begin{equation} \label{eq:BCv2}
    \mathcal{N} \cdot \tau = 
    \begin{bmatrix}
        0 & n_3 & n_2\\
        n_3 & 0 & n_1\\
        n_2 & n_1 & 0
    \end{bmatrix}
    \cdot \tau = 0 \text{ on } \partial \Omega,
\end{equation}
where $\det (\mathcal{N})=2n_1n_2n_3$.  This implies $\tau=0$ at any point on $\partial\Omega$ where the normal is not within a coordinate plane (i.e. $n_i\ne0 \enspace \forall i$).

\begin{figure}
    \centering
    \includegraphics[width=0.6\linewidth]{Omega.jpg}
    \caption{Boundary conditions on the subset $\partial\Omega_1\subset\partial\Omega$ with normal vector orthogonal to the $x_1$-axis (i.e. $n_1=0$).}
    \label{fig:Omega}
\end{figure}

Consider the subset $\partial\Omega_i \subset \partial\Omega$ consisting of all points with $n_i=0$ (see Figure \ref{fig:Omega} for $i=1$).  These points lie in a region of $\partial \Omega$ where derivatives in the $x_i$ direction are tangential to $\partial\Omega$.  The boundary condition \eqref{eq:BCv2} forces $\tau_i=0$ within $\partial\Omega_i$, but the other components are free provided $n_{i+1}\tau_{i+2}+n_{i+2}\tau_{i+1}=0$ where indices are interpreted cyclically.  However, there is further constraint resulting from the form of our general solution as follows.

Since $\tau_i=0$ we have $\partial_iU_{i+1}=\partial_iU_{i+2}$ on $\partial\Omega_i$.  Noting that on $\partial\Omega_i$ the $x_i$ direction is tangent to $\partial \Omega$, we can then compute
\[
\begin{split}
    \partial_i\tau_{i+1}&=\partial_i\partial_{i+1}(-U_{i+2}+U_{i})\\
    &=-\partial_i\partial_{i+1}U_{i+2}\\
    &=-\partial_i\partial_{i+1}U_{i+1}\\
    &=0
\end{split}
\]
on $\partial\Omega_i$, and similarly $\partial_i\tau_{i+2}=0$.  We can carry out this process on any axis-aligned region to conclude that any point on $\partial\Omega$ has either $\tau=0$ or has gradients of $\tau$ equal to zero in at least one tangent direction. It follows from the compact nature of $\Omega$ that we can integrate around $\partial\Omega$ and determine our boundary condition is equivalent to the Dirichlet condition $\tau=0$ on all of $\partial\Omega$. 

Furthermore, we can also consider gradients on $\partial\Omega$ through the decomposition of $\nabla$ into normal and tangential components.  In a local orthonormal basis $\{n,t^1,t^2\}$ we can write the gradient on $\partial\Omega$ in terms of directional derivatives as
\[
\nabla=nD_n+t^1D_{t^1}+t^2D_{t^2}
\]

Given our homogeneous Dirichlet condition, $D_{t^1}\tau=D_{t^2}\tau=0$  on $\partial\Omega$ and hence on the boundary \eqref{eq:ShearEquilibrium} becomes
\[
\mathcal{N}D_n\tau=0.
\]
It follows that $D_n\tau=0$ at any point on $\partial\Omega$ where $n$ is not within a coordinate plane. 

Through a similar process to before, for $n_i=0$ it is possible to show
\[
D_n\tau_i=0 \quad \text{and} \quad \partial_iD_n\tau_{i+1}=\partial_iD_n\tau_{i+2}=0
\]
on $\partial\Omega_i$. Hence, like before, $D_n\tau=0$ on all of $\partial\Omega$ and we conclude that all components of $\tau$ and their derivatives are zero on $\partial\Omega$.  

From \eqref{eq:ChrisGeneral} this means $\partial_{i+1}\partial_{i+2}U_{i}=0$ on $\partial\Omega$ and we can conclude each $U_i$ is additively separated of the form
\[
U_i=V_i(x_{i+1})+W_i(x_{i+2})
\]
within a neighbourhood of the boundary and also throughout $\Omega$ due to the lack of dependence on $x_i$.  From \eqref{eq:ChrisGeneral} the solution must then be of the form $\tau_i(x_i)$ (no Einstein summation implied) and it follows from the Dirichlet condition that $\tau=0$ everywhere within $\Omega$.  

So in summary, if the three diagonal components of $\epsilon$ are zero, it follows that all of the shear components are zero as well.  Hence we have proved our conjecture in the case where we are constructing the full strain tensor from known diagonal components.  

Explicitly; in the context of high-energy x-ray eigenstrain tomography, we have shown that it is sufficient to measure three orthogonal components of normal elastic strain from three orthogonal scalar tomography experiments.  The full elastic strain tensor can be uniquely reconstructed at every point in the sample from this information.

Note that this proof has relied on a number of regularity assumptions relating to $\partial\Omega$ and the fields we have considered.  As mentioned earlier, for completeness, in \ref{apdx:DistProof} we present an additional proof for the more general case where $\partial\Omega$ is only assumed Lipchitz continuous and $\sigma,\epsilon$ are viewed in the weak/distributional sense in the Sobolev space $H^1(\mathcal{S}^2,\Omega)$.  In many ways this more general proof is simpler, however several of the concepts from the first proof are useful for what follows in the rest of the paper.

\section{Minimum requirements of reconstructed eigenstrains}

We now turn our attention to the process of reconstructing stress within $\Omega$ through various numerical schemes focused on determining a corresponding $\epsilon^*$.  This clearly relates to `eigenstrain tomography', but is also common to all inverse eigenstrain problems.  It is important to recognise that this problem is itself ill-posed; multiple eigenstrains can correspond to the same residual stress field \cite{wensrich_well-posedness_2025}.

In practice, the process often relies on implementation of the forward eigenstrain calculation using commercial finite-element software.  The approach is usually to forward map a number of `basis' eigenstrains to corresponding residual stress fields which are then fitted to observations using a least-squares error minimisation.

Often in this setting, eigenstrain is implemented through thermal expansion modelling capability within the software that only has access to diagonal components, and often only isotropically (i.e. hydrostatic eigenstrains).  From this perspective we ask the following pertinent question;

\begin{question}
    What is the minimum number of components of $\epsilon^*$ that need to be considered?  For example, can we span all possible distributions of $\sigma$ with only the diagonal components of $\epsilon^*$ (i.e. by assuming $\epsilon^*_{12}=\epsilon^*_{13}=\epsilon^*_{23}=0$)?  Similarly, what subset of stress distributions can we access if $\epsilon^*$ is assumed to be isotropic at every point?
\end{question}

From this perspective, we state the following theorem;

\begin{theorem}[The span of diagonal eigenstrains]\label{thm:DiagSpan}
On a bounded Lipschitz $\Omega$, say the residual stress field $\sigma$ is generated by the eigenstrain $\epsilon^*$.  Then there exists diagonal eigenstrains $\epsilon^{*d}=\mathrm{diag}(e)$ with $e\in C^2(\mathcal{S}^1,\Omega)$ that map to the same $\sigma$.
\end{theorem}

\begin{proof}
    From \cite{wensrich_well-posedness_2025} we know that the null space of the forward mapping is any eigenstrain of the form
    \[
    \epsilon^{*0}=\nabla_s u
    \]
    where $u\in C^3(\mathcal{S}^1,\Omega)$.  These eigenstrains generate no residual stress.  It follows that our theorem is true if we can find a $u$ such that
    \[
    [\nabla_s u]_{ij} = \epsilon^*_{ij} \text{ for } i<j
    \]
    for any choice of $\epsilon^*$.  In terms of components this becomes
    \begin{align}\label{eq:NullComp}
    \partial_1u_2+\partial_2u_1&=2\epsilon^*_{12}\nonumber\\
    \partial_1u_3+\partial_3u_1&=2\epsilon^*_{13}\\
    \partial_2u_3+\partial_3u_2&=2\epsilon^*_{23}.\nonumber
    \end{align}
    
    Note that the homogeneous solutions to \eqref{eq:NullComp} are, as per \eqref{eq:ChrisGeneral}, of the form
    \begin{equation} \label{eq:HomoPart}
        u^0_i=\partial_i(U_{i+1}-U_{i+2}),
    \end{equation}
    where indices are interpreted cyclically and $U_i$ for each $i$ is an arbitrary function that depends only on two variables with $\partial_iU_i=0$ (no Einstein summation implied). From this perspective, any solution to \eqref{eq:NullComp} will not be unique.
    
    Differentiating \eqref{eq:NullComp} we can write
    \[
    \begin{split}
    \partial_1\partial_3u_2+\partial_2\partial_3u_1&=2\partial_3\epsilon^*_{12}\\
    \partial_1\partial_2u_3+\partial_2\partial_3u_1&=2\partial_2\epsilon^*_{13}\\
    \partial_1\partial_2u_3+\partial_1\partial_3u_2&=2\partial_1\epsilon^*_{23}
    \end{split}
    \]
    within $\Omega$. Rearranging yields
    \[
    \begin{split}
    \partial_2\partial_3u_1&=-\partial_1\epsilon^*_{23}+\partial_2\epsilon^*_{13}+\partial_3\epsilon^*_{12}\\
    \partial_1\partial_3u_2&=\partial_1\epsilon^*_{23}-\partial_2\epsilon^*_{13}+\partial_3\epsilon^*_{12}\\
    \partial_1\partial_2u_3&=\partial_1\epsilon^*_{23}+\partial_2\epsilon^*_{13}-\partial_3\epsilon^*_{12}
    \end{split}
    \]
    which can be integrated directly to determine $u_i$ up to additive fields of the form \eqref{eq:HomoPart}.  We can then compute $e_i=\epsilon^*_{ii}-\partial_iu_i$ (no Einstein summation implied) and hence the theorem is proven.
\end{proof}

From this we form the following Corollary;

\begin{corollary}
    Residual stress fields exist that do not correspond to any isotropic eigenstrain.
\end{corollary}

\begin{proof}
    Consider the following example;
    \begin{equation}\label{eq:Counter}
    \epsilon^*=\begin{bmatrix}
        0 & -\tfrac{1}{2}x_1^2x_3 & 0\\
        -\tfrac{1}{2}x_1^2x_3 & 0 & 0\\
        0 & 0 & 0
    \end{bmatrix}
    \end{equation}
    Following the process above, the corresponding null potentials that link this to a diagonal eigenstrain are of the form
    \[
    u=\begin{bmatrix}
        -\tfrac{1}{2}x_1^2x_2x_3\\
        -\tfrac{1}{6}x_1^3x_3\\
        \tfrac{1}{6}x_1^3x_2
    \end{bmatrix}+
    \begin{bmatrix}
        \partial_1(U_2-U_3)\\
        \partial_2(U_3-U_1)\\
        \partial_3(U_1-U_2)
    \end{bmatrix}
    \]
    resulting in
    \[
    \epsilon^{*d}=\begin{bmatrix}
        x_1x_2x_3-\partial_1^2(U_2-U_3) & 0 & 0\\
        0 & -\partial^2_2(U_3-U_1) & 0\\
        0 & 0& -\partial^2_3(U_1-U_2)
    \end{bmatrix}.
    \]
    Any diagonal eigenstrain that creates the same residual stress as \eqref{eq:Counter} must be of this form.  Given each $U_i$ is only a function of two variables, a choice of $U_i$ that renders all of the diagonal components of this tensor equal does not exist.
\end{proof}

The direct consequence of Theorem \ref{thm:DiagSpan} is that, in any numerical method for inverse eigenstrain problems of any form, we can restrict the search for an eigenstrain corresponding to a set of observations to the three non-zero components of a diagonal tensor (i.e. $e_i$).  This represents a reduction of the dimension of the domain of search by a factor of two.

Furthermore, the corollary implies that it is not possible to reduce the search to one component.  Specifically, this means that approaches that employ isotropic thermal expansion capability in commercial finite element software to solve eigenstrain problems are not robust.

\section{Conclusions and further questions on uniqueness}

Overall we conclude that the `eigenstrain tomography' concept as described by Uzun \emph{et al.} \cite{uzun_tomographic_2024} is well-posed provided that all three diagonal components of elastic strain (or stress) are known/measured. This obviously commands a higher experimental burden than implied by Uzun \emph{et al.}, but only half of the six required tomographic scans implied by Lionheart and Withers \cite{lionheart_diffraction_2015}.  Explicitly; we have shown that it is not possible to reconstruct the full elastic strain (or stress) tensor in three dimensions from measurements of a single component, or even two components.

More broadly, we have also shown that the existence of null eigenstrains (previously termed `impotent' eigenstrains \cite{irschik_eigenstrain_2001,ziegler_eigenstrain_2004}) can be used to our advantage when solving inverse eigenstrain problems of any form.  In particular, we have shown that any eigenstrain has a corresponding set of diagonal eigentrains that generate the same residual stress field.  The direct consequence is that any numerical method focused on finding an eigenstrain corresponding to a set of measurements can, without loss of generality, restrict the search to the three components of a diagonal eigenstrain (rather than six general components).

It is also clear that this process cannot be reduced to finding a single component of hydrostatic (isotropic) eigenstrain.  Many (most) residual stress fields cannot be generated by any isotropic eigenstrain.  Perhaps there are situations where there are grounds to assume the physical eigenstrains within a given sample are of this form, but in general, any numerical method that restricts the search to hydrostatic eigenstrains is not robust.  

Outside of these conclusions, there are a number of other questions that can be posed in this area.  For example; in terms of the x-ray technique discussed in the introduction, each tomographic scan contains much more information that just the normal component of strain along the axis of rotation.  Specifically, each projection consists of the full Transverse Ray Transform of strain in all transverse directions relative to the incident beam.  This raises the question;

\begin{question}
    Is it possible to use all of the transverse strain measurements from high-energy x-ray diffraction along with an assumption of equilibrium and appropriate boundary conditions to uniquely reconstruct $\epsilon$ (hence $\sigma$) from one tomographic scan about a single axis?
\end{question}

Indeed, this was previously carried out numerically by Hendriks \emph{et al.} \cite{hendriks_bayesian_2020} through a machine learning approach, but the question remains: Was this process well-posed?  Were their reconstructions unique?  If they are, this would represent a significant simplification of the measurement process.

More broadly, in terms of more traditional implementation of the eigenstrain method;

\begin{question}
    Often our knowledge involves one or more components of strain or displacement at a number of points or over a number of surfaces (usually planes) within $\Omega$.  In general, how do we characterise the null space of these inverse eigenstrain problems?
\end{question}

And finally;

\begin{question}
    We have shown that, through equilibrium, the full stress tensor is uniquely determined from either the diagonal or off-diagonal components of $\epsilon$.  Is the same true for any set of three distinct components of $\epsilon$?  Is Conjecture \ref{Con:Conj1} true in general?
\end{question}

\bibliography{references}

\appendix

\section{Notation and spaces} \label{Sec:Apndx}

We consider $\Omega\subset\mathbb{R}^3$ as compact with either a $C^1$ smooth or Lipschitz continuous boundary $\partial\Omega$.  In the smooth case (and where it is defined in the Lipschitz case), $n$ is the outward unit normal vector on $\partial\Omega$.

Within $\Omega$ we denote the space of rank-$m$ tensor fields with continuous $l^\mathrm{th}$ derivatives as $C^l(\mathcal{S}^m,\Omega)$. In the case of scalars (i.e. $m=0$), this is denoted $C^l(\Omega)$.  In \ref{apdx:DistProof} we also refer to Sobolev spaces; we denote the Sobolev space of rank-$m$ tensor distributions on $\Omega$ with square-integrable $l^{th}$  weak derivatives $H^l(\mathcal{S}^m,\Omega)$.

Index notation with implied Einstein summation is used throughout the paper except where explicitly stated. Contraction is implied by dots. For example, with $A\in\mathcal{S}^2$ and $b\in\mathcal{S}^1$, $[A\cdot b]_i=A_{ij}b_j$.  Similarly, with $B\in \mathcal{S}^4$, $[B:A]_{ij}=B_{ijkl}A_{kl}$.

We also make use of the following notation and operators;

\begin{itemize}
    \item[--]{The partial derivative with respect to the coordinate $x_i$ is denoted $\partial_i=\frac{\partial}{\partial x_i}$.}
    \item[--]{The divergence operator is $\mathrm{Div}(\cdot)$; for $f\in C^1(\mathcal{S}^2,\Omega)$ we have $[\mathrm{Div}(f)]_j=\partial_i f_{ij}$.}
    \item[--]{The gradient operator is $\nabla$; it is a vector operator with components $\nabla_i=\partial_i$ that is applied as an outer (dyadic) product. For $v\in C^1(\mathcal{S}^1,\Omega)$ we have $[\nabla v]_{ij}=\partial_j v_{i}$.}
    \item[--]{The symmetric gradient operator is $\nabla_s$; for $v\in C^1(\mathcal{S}^1,\Omega)$ we have $[\nabla_s v]_{ij}=\tfrac{1}{2}(\partial_j v_i +\partial_i v_j)$.}
    \item[--]{The scalar directional derivative in the direction of $v$ is $D_v$; for $v\in\mathcal{S}^1$ we have $D_v=\tfrac{v}{||v||}\cdot\nabla$.}
    \item[--]{$\mathrm{diag}(\cdot)$ constructs a diagonal rank-2 tensor from a rank-1 vector; $[\mathrm{diag}(v)]_{ij}=v_i\delta_{ij}$ where $\delta$ is the Kronecker tensor and no Einstein summation is implied}
    
\end{itemize}

\section{`Poorly' vector fields} \label{apdx:poorly}

The reader is likely to be familiar with Killing fields as named after Wilhelm Killing.  These are vector fields $X\in\mathbb{R}^d$ such that $\nabla_s X=0$, or in coordinates
\[
\nabla_iX_j+\nabla_jX_i=0.
\]
In the general case of curvilinear coordinates or on a Riemannian manifold the derivatives are covariant. Killing fields are also called {\em infinitesimal isometries} and in the Euclidean case they define infinitesimal rigid body motions and form a finite dimensional space. 

While Killing fields preserve the metric along with distances and angles, conformal Killing fields preserve only angles and are solutions of
\[
\nabla_iX_j+\nabla_jX_i - \tfrac{2}{d}\nabla_k X_k g_{ij}=0
\]
in dimension $d$ where $g$ is the metric tensor. These also form a finite dimensional space for $d>2$. Consider for comparison volume preserving solenoidal vector fields defined by $\nabla_k X_k=0$, which form an infinite dimensional vector space. 

In this paper we encounter vector fields $X$ (specifically $\tau$ and $u^0$) such that
\[
\nabla_iX_j+\nabla_jX_i=0 \text{ for } i < j,
\]
a weaker condition than both Killing and conformal Killing fields. Not finding these already described in the literature, we name them `poorly fields' on the basis that what does not kill you might still make you poorly. 

We note the following properties;

\begin{itemize}

\item[--]{In contrast to Killing and conformal Killing fields, in the Euclidean case of $d=3$ the space of poorly fields is infinite dimensional; as per \eqref{eq:ChrisGeneral}, they can be represented in terms of three functions of two variables.}

\item[--]{Like Killing and conformal Killing fields, poorly fields form a Lie algebra under the usual Lie bracket of vector fields.}

\item[--]{Killing and conformal Killing fields and volume preserving fields are Lie algebras of infinitesimal automorphisms of $G$-structures, where $G$ is a subgroup of non-singular $d\times d$ matrices \cite{kobayashi_transformation_2012}.  This is not the case for the Lie algebra of poorly fields. The definition of diagonal and off diagonal is coordinate (or at least frame) dependent hence poorly do not have the same kind of invariant properties. }

\item[--]{The famous paper of DeTurke and Yang \cite{deturck_existence_1984} shows that any Reimannian three-manifold has local orthogonal (but not necessarily orthonormal) coordinates. The poorly fields can be viewed as being infinitesimal perturbations of these coordinates that preserve the diagonal structure. The study of coordinate systems in which the metric tensor is diagonal is a problem with a long history; see Zakharov \cite{zakharov_description_1998} and references therein for details.}
\end{itemize}

\section{Proof of Conjecture 1 for diagonal strains in the weak sense.} \label{apdx:DistProof}

\begin{theorem}
Let $\Omega\subset \mathbb{R}^3$ be bounded with $\partial \Omega $ Lipschitz and assume $\sigma\in H^1(\mathcal{S}^2,\Omega)$ is symmetric and satisfies $\mathrm{Div}(\sigma)=0$ inside $\Omega$ and $\sigma\cdot n =0$ on $\partial\Omega$.  If the diagonal components of $\sigma$ are zero then all other components are zero as well.
\end{theorem}
\begin{proof}
We begin by extending $\sigma\in H^1(\mathcal{S}^2,\mathbb{R}^3)$ by defining $\sigma=0$ on $\mathbb{R}^3\backslash\Omega$. Then for any test function $\varphi$ (i.e. smooth and compactly supported in $\mathbb{R}^3$) via integration by parts we have
\[
\begin{split}
\int_{\mathbb{R}^3} \sigma\nabla \varphi dx=&
\int_{\mathbb{R}^3\setminus\Omega} \sigma\nabla \varphi dx + \int_\Omega \sigma\nabla \varphi dx \\
=&\int_{\mathbb{R}^3\setminus\Omega} \sigma\nabla \varphi dx + \int_{\partial\Omega} (\sigma\cdot n) \varphi dS - \int_\Omega \mathrm{Div}(\sigma)\varphi dx =0
\end{split}
\]
where $dx$ is the volume measure on $\mathbb{R}^3$.

Since the diagonal components of $\sigma$ are zero inside $\Omega$ they are also zero in all of $\mathbb{R}^3$.  The remaining off-diagonal components (written in terms of the vector $\tau=[\begin{matrix} \sigma_{23} & \sigma_{13} & \sigma_{12}\end{matrix}]^T$) must then satisfy
\[
\int_{\mathbb{R}^3}(\tau_i\partial_j\varphi+\tau_j\partial_i\varphi)dx=0 \text{ for } i<j
\]
for any choice of $\varphi$. For any test function $\Psi$, the choice of $\varphi=\partial_k\Psi$ is also a test function and hence the above equations also imply
\[
\int_{\mathbb{R}^3} \tau_i\partial_{j}\partial_{k}\Psi dx +\int_{\mathbb{R}^3} \tau_{j}\partial_{i}\partial_{k}\Psi dx = 0
\]
for any $\Psi$ and any distinct choice of $i\ne j\ne k$.  Combining the three resulting equations, it follows that
\begin{equation} \label{eq:transinv2}
\int_{\mathbb{R}^3}\tau_i\partial_{j}\partial_{k}\Psi dx=0.
\end{equation}
for distinct $i\ne j\ne k$.

Let $\varphi$ be an arbitrary test function and let $r>0$ be such that $\varphi(x+re^j)$ is supported outside of $\Omega$ for any index $j$ where $e^j$ is the unit vector aligned with the positive $x_j$ axis.  We now construct a test function as follows;

Let $\varphi_{r_j}(x)=\varphi(x)-\varphi(x+re^j)$ and define
\[
\Psi(x)=\int_{-\infty}^{0}\varphi_{r_j}(x+se^j) ds,
\]
so that $\partial_j\Psi=\varphi_{r_j}$. Note that $\Psi$ is compactly supported.

It follows that
\[
\int_{\mathbb{R}^3}\tau_i\partial_{j}\partial_{k}\Psi dx=\int_{\mathbb{R}^3} \tau_i\partial_k\varphi(x) dx-\int_{\mathbb{R}^3} \tau_i\partial_k\varphi(x+re^j) dx
\]
and since the support of $\varphi(x+re^j)$ is outside the support of $\tau$ we can conclude from \eqref{eq:transinv2} that
\[
\int_{\mathbb{R}^3} \tau_i\partial_k\varphi dx=0.
\]
Through a repeated application of this process (now on the index $k$) it follows that
\[
\int_{\mathbb{R}^3} \tau_i\varphi dx=0
\]
for any test function $\varphi$, and hence in the weak sense the absence of the diagonal components of $\sigma$ imply $\tau=0$ throughout $\mathbb{R}^3$ and specifically within $\Omega$.
\end{proof}

\end{document}